\documentclass{article}
\usepackage{graphicx} % Required for inserting images

\usepackage[T2A]{fontenc}
\usepackage[utf8]{inputenc}
\usepackage[english]{babel}
\usepackage[tbtags]{amsmath}
\usepackage{amsfonts,amssymb,mathrsfs,amscd}

\usepackage{color}

\usepackage{hyperref}

\newlength{\wdth}

\newtheorem{theorem}{Theorem}
\newtheorem{lemma}{Lemma}%[section]
\newtheorem{propos}{Proposition}

\newtheorem{remark}{Remark}

\newif\ifru
\rufalse

\newif\ifen
\entrue % \*true or \*false
\ifen
\title{On strong solution of a multidimensional  SDE: extension of Yamada -- Watanabe's theorem }
\author{A.A. Lyappieva\footnote{M.V. Lomonosov Moscow state university \& Institute for information transmission problems RAS (A.A. Kharkevich Institute), Moscow, Russian Federation; email:  anastasiia.liappieva@math.msu.ru} \, and A.Yu. Veretennikov\footnote{Institute for information transmission problems RAS (A.A. Kharkevich Institute) \& RUDN University, Moscow, Russian Federation; email: alexander.veretennikov2011@ya.ru} }
\fi

\date{}

\begin{document}

\maketitle

\begin{abstract}
A new strong uniqueness result for a multidimensional SDE with a non-degenerate diffusion and partially irregular drift is established. It may be regarded as a combined variation on the  themes of Yamada \& Watanabe (1971), of Zvonkin (1974), and of the second author of the present paper (1980).
\end{abstract}

\section{Introduction}
Consider an SDE in $\mathbb R^d$   ($d\ge 1$)
\begin{equation}\label{SDE1}
dX_t = \sigma(X_t)dW_t + b(X_t)dt, \quad X_0 = x, 
\end{equation}
with a diagonal matrix $\sigma = (\sigma_{ij})$ (that is, with $\sigma_{ij} \equiv 0$ for $j\neq i$), and with dependence $\sigma_{ii}(x^i)$; in addition, each $\sigma_{ii}$ 
is globally 
%? then localization may be mentioned in a remark)} 
uniformly non-degenerated. The drift $b$ is assumed to have a form,  
\begin{equation}\label{b}
b^i(x) = b^i_0(x^i) + b^i_1(x), \quad 1\le i\le d, 
\end{equation}
where $b_1$ is regular in some sense (see what follows), while $b_0$ is only Borel measurable; all coefficients and their parts such as $b_0, b_1$ are assumed bounded. We highlight that $b_1^i$ depends on the whole state $x$ for each $i$, while $b_0^i$ just depends on the component $x^i$, similarly to $\sigma_{ii}$. Of course, the interest is in the case $d>1$, because the case $d=1$ is already completely covered by A.K. Zvonkin's result in \cite{Zvonkin74}. 

\section{Localised Yamada -- Watanabe theorem}
Firstly, let us recall Yamada -- Watanabe's theorem \cite[Theorem 1]{YW71a}, which partial extension will be considered in the present paper. Let us consider an SDE  in $\mathbb R^d$   ($d\ge 1$) more general than (\ref{SDE1}),
\begin{equation}\label{SDE2}
dX_t = \sigma(t,X_t)dW_t + b(t,X_t)dt, \quad X_0 = x, 
\end{equation}
again with a diagonal matrix $\sigma = (\sigma_{ij})$ (that is, with $\sigma_{ij} \equiv 0$ for $j\neq i$), and with a dependence of the coefficients on the time variable: $\sigma_{ii}(t,x^i)$ and $b(t,x)$. All coefficients are bounded and Borel measurable, with certain regularity, which is specified in the following Proposition.

\begin{propos}[\cite{YW71a}]\label{pro1}
Let  $b$ satisfy the inequality 
\begin{equation}\label{rhob}
\sup_{t\ge 0}\sup_{|x-x'|\le u}|b(t,x) - b(t,x')|\le \rho_b(u), \quad \forall\, u\ge 0,
\end{equation}
with the assumption
\begin{equation}\label{irhob}
\int_{0+}\rho_b^{-1}(s)ds = +\infty,
\end{equation}
and also
\begin{equation}\label{rhos}
\sup_{t\ge 0}\max_{i}\sup_{|x^i-(x')^i|\le u}|\sigma_{ii}(t,x^i) - \sigma_{ii}(t,(x')^i)|\le \rho_\sigma(u), \quad \forall\, u\ge 0,
\end{equation}
with the assumption
\begin{equation}\label{irhos}
\int_{0+}\rho_\sigma^{-2}(s)ds = +\infty,
\end{equation}
Then the equation (\ref{SDE2}) has a pathwise unique strong solution.
\end{propos}
{\em NB.} Strictly speaking, the statement of this result in \cite{YW71a} only claims strong uniqueness. However, existence of a weak solution under the assumed conditions follows from Skorokhod's weak existence theorem \cite{Sko} for SDEs with just continuous in $x$ coefficients, and then existence of a strong solution then follows from the {\bf Yamada -- Watanabe principle} \cite[Corollary 3]{YW71a}, see also \cite{KZ}.

\medskip

In the next section, we will need a localized version of this result.  

\begin{propos}[Localized Yamada -- Wanatabe theorem]\label{pro2}
Let for any $R>0$ there exist functions $\rho_{\sigma,R}, \rho_{b,R}$ such that 
\begin{equation}\label{rhobR}
\sup_{t\ge 0}\sup_{|x-x'|\le u, \, |x|, |x'|\le R}|b(t,x) - b(t,x')|\le \rho_{b,R}(u), \quad \forall\, u\ge 0,
\end{equation}
with the assumption for each $R$
\begin{equation}\label{irhobR}
\int_{0+}\rho_{b,R}^{-1}(s)ds = +\infty,
\end{equation}
and also 
\begin{equation}\label{rhosR}
\sup_{t\ge 0}\max_{i}\sup_{|x^i-(x')^i|\le u, \, |x|, |x'|\le R}|\sigma_{ii}(t,x^i) - \sigma_{ii}(t,(x')^i)|\le \rho_{\sigma,R}(u), \quad \forall\, u\ge 0,
\end{equation}
with the assumption that for each $R$
\begin{equation}\label{irhosR}
\int_{0+}\rho_{\sigma,R}^{-2}(s)ds = +\infty.
\end{equation}
Then the equation (\ref{SDE2}) has a pathwise unique strong solution.
\end{propos}

\noindent
{\em Proof. } It follows from proposition \ref{pro1} if we ``smoothly truncate'' the coefficients, e.g., by multiplying them by a Lipschitz function $\chi_R(|x|)$ such that 
$$
\chi_R(u) = 1(|u|\le R) + (R+1 - |u|)\,1(R<|u|\le R+1).
$$
Then we may argue that the SDE (\ref{SDE2}) with new coefficients 
$$
b_R(t,x) = b(t,x)\chi_R(|x|), \; \sigma_R(t,x) = \sigma(t,x)\chi_R(|x|)
$$
has a pathwise unique strong solution, say, $X^R_t, \, t\ge 0$; then due to  pathwise uniqueness for each $R$,  the sequence of these solutions stabilizes as $R\to\infty$; moreover, they all never explode. The limiting process is, clearly, a strong solution of (\ref{SDE2}), which is pathwise unique, as it must coincide with any of $X^R_t$ until the first exit time from the ball $B_R := (x\in \mathbb R^d: \, |x|\le R)$; so it is pathwise unique on the whole half-line $t\ge 0$.
\hfill $\square$

\begin{remark}
It is interesting that a localized version of Ito's theorem is well-known (see \cite[Theorem 4.3.1]{IkedaW}), yet, the authors of the present paper were unable to find a localized version of Yamada \& Watanabe's result. We did it for them now, of course, without pretending to be the authors.
\end{remark}

\begin{remark}
The original proof of proposition \ref{pro1} is based on a special case of the Bihari -- LaSalle lemma. For completeness we provide the proof of this case in the appendix, also without any authorship pretense.

\end{remark}

\section{Main result}
Let us return to the homogeneous equation (\ref{SDE1}). Recall that the drift $b$ has a representation $b = b_0 + b_1$, see (\ref{b}).
\begin{theorem}\label{thm1}
Let $b_1$ satisfy the inequality 
\begin{equation}\label{rhobh}
\sup_{|x-x'|\le u}|b_1(x) - b_1(x')|\le \rho_b(u), \quad \forall\, u\ge 0,
\end{equation}
%so that  
%$$
%\int_{0+}\rho_b^{-1}(s)ds = +\infty,
%$$ 
%(for each $i$, or as a whole function?), 
%in view of Bihari -- LaSalle inequality
with the modulus of continuity function $ \rho_b$ satisfying the assumption (\ref{irhob}), and let $\sigma$ satisfy the inequality
\begin{equation}\label{rhosh}
\max_{i}\sup_{|x^i-(x')^i|\le u}|\sigma_{ii}(x^i) - \sigma_{ii}((x')^i)|\le \rho_\sigma(u), \quad \forall\, u\ge 0,
\end{equation}
with the modulus of continuity function $ \rho_\sigma$ satisfying the assumption (\ref{irhos}).
%$$
%\int_{0+}\rho_\sigma^{-2}(s)ds = +\infty,
%$$ 
In this case, if $\sigma^2$ is uniformly non-degenerate, then the equation (\ref{SDE1}) has a pathwise unique strong solution.
\end{theorem}

\begin{remark}
The ``modulus of continuity'', in particular, means that such a function is continuously increasing on $[0,\infty)$, is concave, and equals zero at the origin. The property of concavity will be explicitly used in what follows.
\end{remark}

\begin{remark}
Recall that due to the result by Yamada and Watanabe \cite{YW71a} (see proposition \ref{pro1} above), if $b_0\equiv 0$, and if all the assumptions of the theorem  \ref{thm1} are met, then the pathwise uniqueness holds, even in the non-homogeneous case where the coefficients may depend on time as well.

Notice that, of course, this result itself may also be localized, but we do not pursue this goal here. In particular, $\sigma$ may be assumed just locally non-degenerate.
\end{remark}

\noindent
{\em Proof.} {\bf 1.} Consider the system of elliptic PDEs, which is, actually, a system  of independent second order ODEs,
\begin{equation}\label{e2}
L^i u^i(x^i) = 0, \; 1\le i\le d, 
%\; \frac{(\sigma_{ii}(t,x^i)}2sigma_{ii}(t,x^i)
\end{equation}
where
$$
L^i = \frac{(\sigma_{ii}(x^i)^2}2 \frac{\partial^2}{(\partial x^i)^2} + b_0^i(x^i)\frac{\partial}{\partial x^i}.
$$
We are interested in a solution $u(x)=(u^1(x^1), \ldots, u^d(x^d))^*$ with all derivative functions $u^i_{x^i}$ strictly positive and locally bounded. This system may be, of course, solved for each $i$ separately and explicitly:
$$
u^i(x^i) = \int_0^{x^i} \exp\left(-2 \int_{0}^{y} \widetilde{b_0^i}(z)dz\right) dy, 
$$
where
$$
\widetilde{b_0^i}(x^i) = \frac{b^i_0(x^i)}{\sigma_{ii}^2(x^i)}.
$$
{\em NB:} We need a solution with the property $(u^i)'>0$, so both integration constants while solving the ODE (\ref{e2}) for each $i$ are chosen zeros, and then, indeed, 
$$
(u^i(x^i) )' = \exp(-2 \int_{0}^{x^i} \widetilde{b_0^i}(y)dy) >0, \quad \forall\, x^i.
$$
Notice that for each $i$, the function 
$$
(u^i(x^i) )'' = -2 \tilde b^i_0 (x^i) \exp(-2 \int_{0}^{x^i} \widetilde{b_0^i}(y)dy) 
$$ 
is locally bounded; the function  $(u^i(x^i) )''/(u^i(x^i) )'$ is bounded globally. 

\medskip

Due to the construction, Ito-Krylov's formula \cite[Chapter 2]{Krylov77} is applicable to $u^i(X^i_t)$, as well as for any other Ito process (i.e., any one with a stochastic differential) with a non-degenerate diffusion coefficient substituted into $u^i(\cdot)$, which process needs not be necessarily a solution of any SDE.

~

{\bf 2.}
We notice that the equation (\ref{SDE1}) does possess a weak solution due to Krylov's result \cite{Krylov69, Krylov77}. We aim to verify pathwise uniqueness. 
So, let us assume that there are two solutions $X_t$ and $Y_t$ with the same Wiener process $B_t$. We are to prove that they coincide. The strategy will be to show it firstly locally, on an interval of time until the exit from a (any) ball of a radius $R$. 

~

Denote
\begin{equation}\label{xxi}
\xi_t:= u(X_t), \quad \eta_t:= u(Y_t), \quad 0\le t\le T.
\end{equation}
By Ito--Krylov's formula, 
\begin{align*}
& d\xi^i_t 
%= du^i(t,X^i_t) 
= (L^i u^i(X^i_t) + b^i_{1}(X_t)u^i_{x^i}(X^i_t))dt + \sigma_{ii}(X^i_t)u^i_{x^i}(X^i_t)dW^i_t 
 \\\\
& =  b^i_{1}(X_t)u^i_{x^i}(X^i_t)dt + \sigma_{ii}(X^i_t)u^i_{x^i}(X^i_t)dW^i_t,
\end{align*}
where the second equality is due to (\ref{e2}).
Similarly, 
\begin{align*}
& d\eta^i_t =  b^i_{1}(Y_t)u^i_{x^i}(Y^i_t)dt + \sigma_{ii}(Y^i_t)u^i_{x^i}(Y^i_t)dW^i_t.
\end{align*}
Notice  that because all derivatives $u^i_{x^i}$ are positive, the mapping $x\mapsto u(x)$ is invertible, so that 
there exists an inverse function 
$v(\xi):= u^{-1}(\xi)$. 
Hence, we can rewrite both SDEs for $\xi_t$ and $\eta_t$ without references to $X_t$ and $Y_t$:
\begin{align}\label{exi}
& d\xi^i_t = \hat b^i(\xi_t) dt + \hat \sigma_{ii}(\xi^i_t) dW^i_t, \end{align}
and
\begin{align}\label{eeta}
& d\eta^i_t =  \hat b^i (\eta_t)dt + \hat \sigma_{ii}(\eta^i_t)dW^i_t.
\end{align}
Here
\begin{align*}
 \hat b^i(y) = b^i_{1}(v(y))u^i_{x^i}(v^i(y^i)), 
\end{align*}
and
\begin{align*}
\hat \sigma_{ii}(y^i) = \sigma_{ii}(v^i(y^i))u^i_{x^i}(v^i(y^i)).
\end{align*}
%{\color{blue}
%So, 
%\begin{align*}
%& d(\xi^i_t - \eta^i_t) =  (b^i_{1}(X_t)u^i_{x^i}(X^i_t) - b^i_{1}(Y_t)u^i_{x^i}(Y^i_t))dt 
% \\\\
%&+ (\sigma_{ii}(X^i_t)u^i_{x^i}(X^i_t) - \sigma_{ii}(Y^i_t)u^i_{x^i}(Y^i_t)) dW^i_t
% \\\\
%&= [\hat b^i(\xi_t) - \hat b^i(\eta_t)]dt + [\hat \sigma_{ii}(\xi^i_t) - \hat \sigma_{ii}(\eta^i_t)] dW^i_t.
%\end{align*}
%}
Let us evaluate the moduli of continuity of $\hat b$ and $\hat \sigma$ on compacts. 
Let $R > 0$, $B_R = \{ y \in \mathbb{R}^d : |y| \leq R \}$, and $K_R = v(B_R)$.

On the compact $K_R$, the functions in the definitions of $\hat{b}$ and $\hat{\sigma}$ possess the following properties:
\begin{enumerate}
\item 
$
%b_1^i (y) = 
\hat b^i(y) = b^i_{1}(v(y))u^i_{x^i}(v^i(y^i))$ is bounded and has a modulus of continuity 
$$
\rho_{\hat b^i,R}(u) \le C_R\rho_b(u) + \|b_1\|_B L^{u^i_{x^i} \circ v^i}_R u,
$$
with $C_R=\sup_{|y|\le R} |u^i_{x^i}(v^i(y^i))|\,L^v_R$, where $L^v_R$ is the Lipschitz constant for the function $v$ on $B_R$, and $L^{u^i_{x^i} \circ v^i}_R$ is the Lipschitz constant for the function $u^i_{x^i} \circ v^i$ on $B_R$. Thus, 
\begin{equation}\label{rhobhat}
\rho_{\hat b,R}(u) \le C_R\rho_b(u) + \|b_1\|_B \max_i L^{u^i_{x^i} \circ v^i}_R u =: C_R\rho_b(u) + \|b_1\|_B L^{u_{x} \circ v}_R u.
\end{equation}

\item
Similarly, $\hat \sigma_{ii}$ possesses a modulus of continuity satisfying
$$
\rho_{\hat \sigma_{ii},R}(u) \le C_R\rho_b(u) + \|b_1\|_B L^{u^i_{x^i} \circ v^i}_R u,
$$
with the same constant $C_R=\sup_{|y|\le R} |u^i_{x^i}(v^i(y^i))|\,L^v_R$, where $L^v_R$ is the same Lipschitz constant for the function $v$ on $B_R$, and $L^{u^i_{x^i} \circ v^i}_R$ is the Lipschitz constant for the function $u^i_{x^i} \circ v^i$ on $B_R$. Thus, 
$$
\rho_{\hat \sigma(u),R} \le C_R\rho_b(u) + \|\sigma\|_B \max_i L^{u^i_{x^i} \circ v^i}_R u =: C_R\rho_\sigma(u) + \|b_1\|_B L^{u_{x} \circ v}_R u.
$$
\end{enumerate}

~

{\bf 4.}
Let us verify that the coefficients $\hat b$ and $\hat \sigma$ satisfy the assumptions of proposition \ref{pro2}\footnote{But, in general, not the assumptions of proposition \ref{pro1}, because of a local nature of the continuity properties of $\hat b$ and $\hat \sigma$}.
We have to check two equalities: 
\begin{equation}\label{divb}
\int_{0+} \frac{1}{\hat{\rho}_{\hat b,R}(s)} ds = +\infty, 
\end{equation}
and
\begin{equation}\label{divs}
\int_{0+} \frac{1}{(\hat{\rho}_{\hat\sigma,R}(s))^2} ds = +\infty.
\end{equation}
Both of them follow from the conditions on $\rho_b$ and $\rho_\sigma$, respectively, and from the concavity of any modulus of continuity. Indeed, from (\ref{rhobhat}) we obtain that 
there are only two options: either $C_R\rho_{b}(u)\le Cu, \, \forall u\ge 0$ with $C=\|b_1\|_B L^{u_{x} \circ v}_R$, or there exists $u_0>0$ such that  $C_R\rho_{b}(u_0)> Cu_0$, and in this case $C_R\rho_{b}(u)\ge Cu, \, \forall 0\le u\le u_0$. In both cases the integral under consideration diverges. 
In the first case it holds 
since 
$$
\int_0^1 \frac{du}{C_R\rho_{b}(u)+Cu} \ge 
\int_0^1 \frac{du}{2Cu} =\infty,
$$
while in the second case it is true due to the assumption (\ref{irhob}) and because 
$$
\int_0^1 \frac{du}{C_R\rho_{b}(u)+Cu} \ge 
\int_0^{1\wedge u_0} \frac{du}{2C_R\rho_b(u)} =\infty.
$$
This shows (\ref{divb}).
Quite similarly (\ref{divs}) is justified, with the only change that here the result is either 
$$
\int_0^1 \frac{du}{(C_R\rho_{b}(u)+Cu)^2} \ge 
\int_0^1 \frac{du}{(2Cu)^2} =\infty,
$$
or (with a new $u_0>0$)
$$
\int_0^1 \frac{du}{(C_R\rho_{\sigma}(u)+Cu)^2} \ge 
\int_0^{1\wedge u_0} \frac{du}{(2C_R\rho_{\sigma}(u))^2} =\infty.
$$
The theorem is proved. \hfill $\square$

\section*{Appendix}
Here is the  ``zero case'', or the simplified version of the Bihari -- LaSalle inequality (see \cite{Bihari,Bihari2})  with a short proof, on which the proof of uniqueness is based. It is needless to say that we do not claim authorship in any way. The purpose of this appendix is just to help the reader, because we were unable ourselves to find a reference with such an explicit statement and with a proof.

\begin{lemma}[simplified Bihari -- LaSalle lemma]\label{Lem-bihari}
If $(v(x), \, x\ge 0)$ is a finite continuous function satisfying 
\begin{equation}\label{a-bihari}
0\le v(x) \le \int_0^x \phi(v(s))ds, \quad x\ge 0, 
\end{equation}
with a nonnegative nonstrictly increasing function $\phi$ such that $\phi(0)=0$ and $\displaystyle \int_0^1 \frac{ds}{\phi(s)} =\lim_{a\downarrow 0} \int_a^1 \frac{ds}{\phi(s)}=\infty$, then 
\begin{equation}\label{v0}
v(x) \equiv 0.
\end{equation} 
\end{lemma}
\noindent
{\em Proof.} 
Denote
\begin{equation}\label{F}
F(x):= \int_0^x \phi(v(s))ds, \quad x\ge 0, 
\end{equation}
so that the condition of the lemma can be rewritten as 
\begin{equation}\label{vF}
0\le v(x) \le F(x). 
\end{equation}
Since $\phi$ is increasing, it also follows that
\begin{equation}\label{viF}
0\le v(x) \le \int_0^x \phi(F(s))ds.
\end{equation}
From the definition of $F$ and from the assumptions it follows that $F$ is increasing as well. 

Now, if $\phi(F(x)) = 0$ for any $x\ge 0$, then clearly $v(x)=0$ from the condition and from the last inequality; hence, (\ref{v0}) is valid. Assume that $\phi(F(x))>0$ for all values  $x>x^*$ for some $x^*\ge 0$, and $\phi(F(x^*))=0$; naturally, in this case $\phi(F(x))=0, \, \forall \, 0\le x\le x^*$. By virtue of (\ref{F}), (\ref{vF}),  and (\ref{viF}), it then follows that $F(x^*)=0$.

Next, notice that $F$ is differentiable and 
$$
0\le F'(x) = \phi(v(x)) \stackrel{\text{by the condition}}\le\phi(F(x)), \quad \forall \, x\ge 0.
$$
Then for any $x>x^*$ we have
$$
0\le \frac{F'(x)}{\phi(F(x))} \le 1.
$$
So, for any $\epsilon>0$ and $x\ge x^* +\epsilon$ we estimate
$$
0\le \int_{x^*+\epsilon}^x \frac{F'(s)ds}{\phi(F(s))} \le x-x^*-\epsilon.
$$
Denoting $F(s)=:z$ and changing variables, the latter inequality may be presented as
$$
0\le \int_{F^{}(x^*+\epsilon)}^{F(x)} \frac{dz}{\phi(z)} \le x-(x^*+\epsilon).
$$
However, as we have already seen, $F(x^*)=0$, and so, $F^{}(x^*+\epsilon) \to 0$ as $\epsilon \to 0$. Therefore, by virtue of the condition $\displaystyle \int_0^1 ds/\phi(s) =\infty$, we get in the limit 
$$
\infty \le x-x^*, 
$$
which is not possible since both $x^*$ and $x$ are finite. This shows that the assumption that $\phi(F(x))>0$ for $x>x^*$ with some $x^*\ge 0$ is not possible either, which means that $\phi(F(x))=0 \;\, \forall \, x\ge 0$. We have already seen that in this case 
$$
v(x) = 0,\quad x\ge 0, 
$$ 
due to the assumption (\ref{a-bihari}), as required\footnote{The ``full'' Bihari -- LaSalle inequality states some bound under the assumption that $v(0) \ge 0$, and the usual form of this inequality is stated for $v(0)>0$. Note that in the lemma \ref{Lem-bihari} it follows automatically from the assumption that $v(0)=0$. {\color{blue}We will happily remove this appendix if anybody, a referee, or someone else, advises us a correct reference to an explicit statement of this precise lemma with a proof, not to the general Bihari -- LaSalle inequality where the ``zero case'' is not directly specified.}}. \hfill $\square$

\end{document}